\documentclass[12pt]{amsart} 
\usepackage{amssymb}

\usepackage{helvet} 

\usepackage{courier} 

\reversemarginpar 
 \normalmarginpar 

\let\oldmarginpar\marginpar 
\renewcommand\marginpar[1]{\-\oldmarginpar{\raggedright\small\sf #1}}

\title[Chow groups and higher congruences]{Chow groups and higher congruences 
for the number of rational points on proper varieties over finite fields}

\author{Najmuddin Fakhruddin}
 
\address{School of Mathematics, Tata Institue of Fundamental Research,
Homi Bhabha Road, Mumbai 400005, India}

\email{naf@math.tifr.res.in}

\newcommand{\nc}{\newcommand}
\nc{\rc}{\renewcommand}

\newcommand{\bs}{\backslash} 
\newcommand{\te}{\otimes}
\newcommand{\bP}{\mathbb{P}} 
\newcommand{\Q}{\mathbb{Q}}
\newcommand{\Z}{\mathbb{Z}}

\newcommand{\A}{\mathbb{A}}
 
\nc{\W}{\mathbb{W}}
\newcommand{\mc}{\mathcal} 
\newcommand{\mb}{\mathbb}

\nc{\ql}{\mathbb{Q}_{\ell}}

\nc{\G}{\Gamma}

\renewcommand{\P}{\mb{P}}

\nc{\aff}{{\A}^1} 
\nc{\naive}{\!\sim_n} 
\nc{\Spec}{\mathrm{Spec}}

\nc{\wt}{\widetilde} 
\nc{\wh}{\widehat}

\nc{\sr}{\stackrel}

\nc{\lr}{\longrightarrow} 
\nc{\wb}{\overline}

\nc{\rank}{\operatorname{rank}}
\nc{\chr}{\operatorname{char}}

\newtheorem{thm}{Theorem}[section] 
\newtheorem{prop}[thm]{Proposition}
 
\newtheorem{cor}[thm]{Corollary}
\newtheorem{lem}[thm]{Lemma}

\theoremstyle{definition}

\newtheorem{rem}[thm]{Remark}
 
\newtheorem{ex}[thm]{Example}

\input{xypic}

\begin{document}

\begin{abstract}
Given a proper family of varieties over a smooth base, with 
smooth total space and general fibre, all over a 
finite field $k$ with $q$ elements, we show that a finiteness hypotheses on the
 Chow groups, $CH_i$,  $i=0,1, \dots , r$,
 of the fibres in the family leads to congruences $\mod q^{r+1}$ for the number
of rational points in \emph{all} the fibres over $k$-rational
points of the base. These hypotheses on the Chow groups are
expected to hold for families of low degree intersections in
many Fano varieties leading to a broad generalisation of the
theorem of Ax--Katz \cite{katz-ax}, as well as results of the
author and C.~S.~Rajan \cite{nf-rajan}. As an unconditional
application, we give an asymptotic generalisation of the Ax--Katz
theorem to low degree intersections in a large class of homogenous spaces.

\end{abstract}

\maketitle

\section{Introduction}
Let $X$ be a smooth proper geometrically connected
variety over a finite field $k$ with $q$ elements and
assume that the base change maps of \emph{rational} Chow groups,
$CH_i(X) \to CH_i(X_{\wb{k(X)}}) $
are isomorphisms for $i = 0,1,\dots,r$ and some $r \geq 0$.
These isomorphisms give rise to a
decomposition of the diagonal correspondence in $X \times X$
(see Lemma \ref{lem:decomp} for details),
which along with the Grothendieck--Lefschetz trace formula shows that
$|X(k)| \equiv \sum_{i=0}^r b_{2i}q^i \mod q^{r+1}$, where
 $b_i = \dim H^i_{et}(X_{\bar{k}}, \ql)$ with $\ell$ a prime not equal to $\chr(k)$;
 this was first observed
by H.~Esnault \cite{esnault-fano}. 
The goal of this paper is to prove
similar results in the case where $X$ is not necessarily smooth: that some such
results should hold under fairly general hypotheses
is suggested by the Ax--Katz theorem  \cite{katz-ax} which gives congruences
as above for intersections of  low degree hypersurfaces in $\bP^n$. For $r=0$,
such results were  proved by the author and C.~S.~Rajan 
in \cite{nf-rajan} under very general hypotheses; 
here we will show that a refinement
of those methods implies congruences for singular varieties also when $r>0$. 
Our results can thus be viewed
as a motivic explanation for, as well as a broad generalisation of (albeit with some
extra hypotheses which should conjecturally always hold), the theorem of Ax and
Katz.

For singular varieties, even for $r=0$, the Chow groups do not
``control'' congruences. The main observation of \cite{nf-rajan} was that if we consider
a family of varieties over a smooth base with smooth total space, then $CH_0$ 
of the generic
fibre (over a sufficiently large extension of the base field) 
does in fact ``control'' congruences $\mod q$ for the number of rational points
of all closed fibres. If
$r>0$ it is easily seen, for example by blowing up a
subvariety of the total space not intersecting the generic fibre,
that a condition solely in terms of the Chow groups of the
generic fibre cannot suffice to imply congruences $\mod q^{r+1}$ for closed
fibres; some extra geometric and/or cohomological conditions must be added.
The main result of this paper is then the following:

\begin{thm}
Let $f: X \to Y$ be a proper dominant morphism of smooth irreducible
varieties  with smooth geometrically irreducible
generic fibre $Z$ over a  finite field $k$ with $q$ elements.
For each point $y$ of $Y$, let $Z_y = f^{-1}(y)$ and
let $\mc{Y} = \Spec(\mc{O}_{y,Y})$.
Assume that for each $i = 0, 1, \dots, r$, the specialisation maps
(\cite[p.399]{fulton-it})
\[
CH_{i}(X \times_Y (\mc{Y} - \{y\})/ \mc{Y}- \{y\}) 
 \to CH_i(Z_y) 
\]
are surjective and the maps
$CH_i(Z_y)   \to CH_i((Z_y)_K)$
are isomorphisms  
for some embedding $k(y) \hookrightarrow K$ with $K$ a universal domain and
$i = 0,1,\dots,r$.
Further, assume that $R^jf_* \ql$ is a constant local system on $Y$ for 
$j = 0,1, \dots, 2r$. Then for all $y
\in Y(k) $, 
\[
|Z_y(k)| \equiv \sum_{i=0}^{r} b_{2i} \ q^i \mod q^{r+1} \ , 
\]
where $b_{2i} = \rank R^{2i}f_* \ql$. \label{thm:main}
\end{thm}

Some condition similar to ours on the $R^jf_* \ql$ probably follows from 
the conclusion of the theorem (for all finite extensions of $k$); it appears to be
close to the minimal extra condition not pertaining to the generic fibre needed
to obtain congruences mod $q^{r+1}$ for $r>0$. But it seems likely that
the hypotheses on the Chow groups of non-generic fibres is implied
by the above cohomological condition along with the hypothesis on 
the Chow groups of the generic fibre.

\smallskip
A theorem of K.~Paranjape \cite{paranjape}
(which was improved in \cite{ELV}) 
implies that the hypotheses of the theorem on
Chow groups (for a fixed $r$) hold
for the universal family of intersections of 
hypersurfaces of multidegree $(d_1,d_2, \dots, d_t)$
in $\bP^n$ if $n$ is sufficiently large compared to $\sum_{k=1}^t d_k$; the other
conditions are easily checked. Thus we recover an asymptotic version of the
Ax--Katz theorem and the full Ax--Katz theorem would
follow from standard conjectures relating the
coniveau filtration on cohomology and Chow groups.
By partly extending Paranjape's results to homogenous spaces we are able to obtain the
following unconditional result:
\begin{thm}
\label{thm:assaxkatz}
Let $k$ be  finite field with $q$ elements, let
$G_n$ be a $k$-form of one of the simple algebraic groups $SL_{n+1}$, $Sp_{2n}$, $SO_{2n}$
or $SO_{2n+1}$ and let
$P_n$ be a maximal parabolic subgroup of $G_n$ also defined over $k$ except for
the case where $G_n$ is $Sp_{2n}$ and $P_n$ is the stabilizer of a maxmimal isotropic
subspace. 
Let $X_n = G_n/P_n$ and $L_n$ the ample generator of $Pic(X_n)$. 
For any tuple  $D := (d_1,d_2, \ldots, d_t)$  of positive integers and any
positive integer $r$, there exists an integer $N(D,r)$ such that for 
any $n > N(D,r)$ and any $s_i \in H^0(X_n, L_n^{\otimes d_i})$, $i=1,2, \ldots,t$
\[
|Z(s_1,s_2,\ldots,s_t)(k)| \equiv \sum_{j=0}^{r-1}\dim(b_2(X_n))q^j \mod q^r
\]
where $Z(s_1,s_2,\ldots,s_t)$ is the subscheme of $X_n$ given by the vanishing of all
the $s_i$'s.
\end{thm}

We also include an example (\ref{ex:E8}) of a homogenous space for the exceptional group $E_8$
for which one also has ``higher conguences'' for the number of rational points
on certain hypersurfaces.

\bigskip
The main idea of the proof of Theorem \ref{thm:main}
is to use the hypothesis on the
Chow groups to obtain a suitable decomposition of the diagonal correspondence
in $X \times_Y X$ which, using standard properties of etale cohomology,
implies the divisibility of the eigenvalues of Frobenius acting on a
certain quotient of $H^*_c(X, \ql)$. An application of the Grothendieck--Lefschetz
trace formula completes the proof. The decomposition is obtained by applying
Lemma \ref{lem:decomp} to the generic fibre which gives us the desired 
decomposition except for an ``error term'' which is a 
correspondence whose support
does not intersect the generic fibre. This correspondence is controlled by applying
Lemma \ref{lem:decomp} inductively  to certain non-generic fibres 
along with a ``nilpotence'' argument  which
eventually shows that the error term can be neglected.

\bigskip
\emph{Notations and Conventions.} By a variety we shall mean a reduced separated
scheme of finite
type over a field $k$. For $X$ a variety, we shall abuse notation and denote by
$H^*(X, \ql)$ the etale cohomology groups $H^*_{et}(X_{\bar{k}}, \ql)$ where $\bar{k}$ is
an algebraic closure of $k$, and similarly for $H^*_c(X, \ql)$, the 
etale cohomology groups with compact supports;
these are finite dimensional $\ql$-vector spaces with a continuous action
of the Galois group of $\bar{k}$ over $k$. All Chow groups
will be considered with rational coefficients.

\section{Lemmas} 

\label{sec:lemmas}

In this section we prove several lemmas which will be used to
prove our main theorem in the next section.

\begin{lem}
\label{lem:decomp}
Let $W$ and $Z$ be varieties over a field $k$ and let $K$ be a universal
domain over $k$. Suppose the maps 
$CH_i(Z)  \to CH_i(Z_K) $
are isomorphisms for $i=0,1, \dots, r$.
Then any $\alpha \in CH_d(W \times Z) $ can be written as
\[
\alpha = c \sum_{i=0}^r \sum_{j=0}^{m_i} [A_{i,j}\times B_{i,j}]  + \beta
\]
where $A_{i,j}$, $B_{i,j}$ are irreducible subvarieties of $W$ and $Z$
respectively of respective dimensions
$d-i$ and $i$,  $\beta$ is a cycle on $W \times Z$ such that
$p_1(|\beta|)$ has dimension $\leq d -r -1$ and $c \in \Q$.
\end{lem}

\noindent Note that we do not assume that $W$ and $Z$ are smooth or even irreducible.

\begin{proof}
This follows using standard methods (originally due to S.~Bloch for $r=0$)
but we give a  proof for the 
reader's convenience. 

The hypothesis implies that for any extension $k'$ of  $k$, the maps
$CH_i(Z) \to CH_i(Z_{k'})$ is an isomorphism in the given range; in 
particular this holds for the function fields, say  $k_1, k_2, \dots, k_m$,
of the irreducible components of $W$. Let 
$\iota_l: k_l \times_k Z \to W \times Z$ be the maps induced by the inclusion of
the generic points and let $\alpha_l = \iota_l^*(\alpha) \in CH_{d - d_l}$,
where $d_l$ is the dimension of the irreducible component of $W$ corresponding
to $k_l$ and  $l =1,2 \dots, m$. By hypothesis,
if $d - d_l \leq r$, then $\alpha_l$ is in the image of the map 
$CH_{d -d_l}(Z) \to CH_{d-d_l}(Z_{k_l})$ (say $\wt{\alpha_l} \mapsto \alpha_l$),
so
by taking  the sum of the Zariski closures of (representatives for) 
the cycles $k_l \times_k \wt{\alpha_l}$  for $l$ such that $d - d_l \leq r$ 
we get a cycle $\gamma$ on $W \times Z$ which is a linear combinations of terms of
the form  $[A \times B]$ with $\dim(A) + \dim(B) = d$,
$\dim(B) \leq r$ and with $A$  one of the irreducible components of $W$.
By construction, $\gamma$ has the property that
$\alpha = \alpha' + \gamma + \beta'$ where  $\alpha'$  is supported
on $W' \times Z$ with $W'$ is a closed subvariety of $W$ such that
 $W \bs W'$  Zariski dense in $W$ 
and $\beta'$ has the same support properties
as those desired of $\beta$. 
We may then replace $W$ by $W'$ and $\alpha$ by $\alpha'$ and complete
the proof by descending Noetherian induction.
\end{proof} 

For the remainder of this section we assume that 
the reader is familiar with the methods of \cite{nf-rajan}, in particular the
notion of proper correspondences and their action on \'etale cohomology
from \cite[Section 2]{nf-rajan}.

\begin{lem}
\label{lem:compos}
Let $f: X \to Y$ be a proper dominant morphism of smooth irreducible
varieties over a field $k$ with smooth geometrically irreducible
generic fibre $Z$. Let $A$, $B$ be subvarieties of $X$ mapping dominantly
to $Y$.
Then the image of
the proper correspondence $[A \times_Y B]$ acting on $H^*_c(X, \ql)$
is contained in the image of 
\[
H^*_c(Y, \ql) \otimes [A] \subset H^*_c(Y, \ql) \otimes H^*(X, \ql) 
\sr{f^* \otimes Id}{\lr}  H^*_c(X, \ql) \otimes H^*(X, \ql) \sr{\cup}{\lr} H^*_c(X, \ql) \ .
\]
\end{lem}

\begin{proof}
Using the inclusion of  $A$ in $X$ and the morphism $f$ we may view $A$
 as giving a proper correspondence in $X \times Y$ and
similarly we may view $B$ as a proper correspondence in $Y \times X$.
 $[A \times_Y B]$ is then
the composite of these two correspondences. It is clear 
from the definitions that the image
of the first correspondence (mapping $H^*_c(Y, \ql) $ to $H^*_c(X, \ql)$ )
is exactly the image of 
$H^*_c(Y, \ql) \otimes [A]$ in $H^*_c(X, \ql)$. The lemma then follows 
from the compatibility of the composition of correspondences with the
action of correspondences on cohomology.
\end{proof}

\begin{lem}
\label{lem:moving}
Let $f: X \to Y$ be a proper dominant morphism of smooth irreducible
varieties over a field $k$ with smooth geometrically irreducible
generic fibre $Z$. For each point $y$ of $Y$, let $Z_y = f^{-1}(y)$ and
let $\mc{Y} = \Spec(\mc{O}_{y,Y})$.
Assume that for each $i = 0, 1, \dots, r$, the specialisation maps
(\cite[p.399]{fulton-it})
\[
CH_{i}(X \times_Y (\mc{Y} - \{y\})/ \mc{Y}- \{y\}) 
 \to CH_i(Z_y) 
\]
are surjective. Let $Y'$ be an
irreducible proper subvariety of $Y$ and let $A_j, B_j$, $j =1,2$, be irreducible
subvarieties of $X$ mapping dominantly to $Y'$. 
Suppose further that $s = \dim(B_1) - \dim(Y') \leq r$.
Then the  proper correspondence $[A_1 \times_Y B_1] \circ [A_2 \times_Y B_2]$
on $X \times_Y X \subset X \times X$ is equivalent (in the group of proper
correspondences) to a correspondence $\Gamma $ supported in
$X \times_Y X$ and such that the image of $|\Gamma|$ in $Y$
is contained in a proper closed
subvariety of $Y'$.
\end{lem}

\begin{proof}
As in the proof of Lemma \ref{lem:compos} we may view the $A_j$'s as
proper correspondences in $X \times Y$ and the $B_j$'s as proper
correspondences in $Y \times X$. It then suffices to show that
$[B_1] \circ [A_2]$, which is a priori a correspondence in $Y \times Y$
supported in the diagonal, is in fact supported on a proper closed
subset of $Y'$ (diagonally embedded in $Y \times Y$).

Let $y$ be the generic point of $Y'$. Then $B_1$ gives an element in 
$CH_s(Z_{y})$ which, using the surjectivity of the specialisation maps,
must be the specialisation of an element 
$ \alpha \in CH_s(X \times_Y (\mc{Y} - \{y\})/ \mc{Y} - \{y'\})$ .
This implies that the element which $B_1$ gives in
$CH_{s}( X \times_Y \mc{Y} )/\mc {Y})$ is $0$:
this is clear if $\mc{O}_{y,Y}$ is a d.v.r.~and we can reduce to
this case  (since $\dim(\mc{O}_{y,Y}) \geq 1$) using functoriality of 
the Gysin maps of \cite{fulton-it}
for a composite of regular embeddings. 
Hence $B_1$ is rationally equivalent, in $X$, to a cycle
$\beta$ such that $f(|\beta|) \cap Y'$ is a proper closed subset
of $Y'$. Then $\beta$ may be viewed as  a proper correspondece in 
$Y \times X$ equivalent to $[B_1]$. This then implies the claim
on the support of $[B_1] \circ [A_2]$.
\end{proof}

\begin{lem}
\label{lem:inj}
Let $f: X \to Y$ be a proper dominant morphism of smooth irreducible
varieties over a  field $k$ with smooth geometrically irreducible
generic fibre $Z$. 
Let $\alpha_{i,k} \in H^{2i}(X, \ql)$, $i =0,1, \dots, r$, $k = 1, \dots, b_{2i}$,
be algebraic cohomology classes such that their restrictions
 to the generic fibre form a basis for
$\oplus_{i=0}^{r} H^{2i}(Z, \ql)$. Then the map
\[
\big( \bigoplus_{i,k} H^*_c(Y, \ql) \big ) \sr{\oplus f^* \cup \alpha_{i,k}}{\lr} H^*_c(X, \ql)
\]
is injective.
\end{lem}

\begin{proof}
Since the $\alpha_{i,k}$'s are algebraic classes and $Z$ is smooth and proper,
we may find algebraic 
cohomology classes
$\beta_{i,k} \in H^*(X, \ql)$ such that $\alpha_{i,k} \cup \beta_{i',k'}$
restricts to $\delta_{i,i'} \cdot \delta_{k,k'} \cdot [pt] \in H^{2\dim(Z)}(Z, \ql)$,
where $\delta$ is the Kronecker delta and $[pt]$ means the cohomology class 
of a point: this follows immediately from  the Hard Lefschetz
theorem if $A$ is projective and one reduces the general case to this by using Chow's lemma,
alterations and the projection formula.

It follows from the construction that  $f_*( \alpha_{i,k} \cup \beta_{i',k'}) = 
\delta_{i,i'} \cdot \delta_{k,k'} \cdot [Y] \in H^0(Y, \ql)$.
For any cohomology class $\gamma \in H^*_c(Y, \ql)$,  it then follows 
from the projection formula that
$f_*(f^*(\gamma) \cup \alpha_{i,k} \cup \beta_{i',k'}) = \delta_{i,i'} \cdot \delta_{k,k'} \cdot
\gamma$, which proves the injectivity.

\end{proof}

\begin{lem}
\label{lem:div}
Let $X_1, X_2$ be smooth geometrically irreducible varieties over a finite
field $k$ with $q$ elements. Let $\Gamma$ be a proper correspondence in
$X_1 \times X_2$ such that $\dim(p_1(|\Gamma|) < \dim(X) - r$. Then all
the eigenvalues of the geometric Frobenius acting on 
$\Gamma^*(H^*_c(X_2, \ql)) \subset H^*_c(X_1, \ql)$ are 
algebraic integers divisible by $q^{r+1}$.
\end{lem}
\begin{proof}
This is a direct generalisation of \cite[Proposition 2.4]{nf-rajan} so we do not give
a detailed proof. The key point is that if $\pi: W_1 \to W_2$ is a map of smooth
geometrically irreducible varieties over $k$ such that $\dim((W_1) < \dim(W_2) - r$
then all the eigenvalues of geometric Frobenius acting on
$\pi_*(H^*_c(W_1, \ql)) \subset H^*_c(W_2, \ql)$ are algebraic integers 
divisible by $q^{r+1}$.
\end{proof}

\section{Proof of the main theorem}

\begin{proof}[Proof of Theorem \ref{thm:main}]
We apply Lemma \ref{lem:decomp} with $\alpha = [\Delta_Z]$,
the class of the diagonal in $Z \times Z$. Let $\wb{A}_{i,j}, \wb{B}_{i,j}$ denote the
Zariski closures of $A_{i,j}, B_{i,j}$ in $X$ and let $\wb{\beta}$ be the Zariski closure of the
cycle $\beta$ in $X \times_Y X$. It follows that $\Delta_X$, viewed as a proper correspondence
in $X \times X$, differs from $c \sum_{i=0}^r \sum_{j=0}^{m_i}[\wb{A}_{i,j} \times_Y \wb{B}_{i,j}]
+ \wb{\beta}$ by a proper correspondence $\Gamma$ supported in $X \times_Y X$ and such that
$p_1(|\Gamma|) = p_2(|\Gamma|)$ does not contain the generic point of $Y$. Applying
Lemma \ref{lem:decomp} inductively to the fibres over the generic points of 
$p_1(|\Gamma|)$ and the induced cycles on these fibres we see that we may write
\[
[\Delta_X] = c \sum_{i=0}^r \sum_{j=0}^{m_i}[\wb{A}_{i,j} \times_Y \wb{B}_{i,j}]
+ d \sum_{s=0}^r \sum_{t=0}^{m_s} [C_{s,t} \times_Y D_{s,t}] + \Gamma'
\]
where the $C_{s,t}, D_{s,t}$ are subvarieties of $X$ with $f(C_{s,t}) = f(D_{s,t}) \subsetneqq Y$
and $\dim(D_{s,t}/f(D_{s,t})) = s$. Furthermore, $\Gamma'$ is a proper correspondence
supported in $X \times_Y X$ with $\dim(p_1(|\Gamma'|)) < \dim(X) -r$.

Now repeatedly  using the fact that $[\Delta_x]^2 = [\Delta_X]$ as proper correspondences and
Lemma \ref{lem:moving}, we see that $[\Delta_X]$ can be written as a sum of
products of proper correspondences, all of them supported in $X \times_Y X$, such that
each product contains a factor of one of the following two types:
(1) $[\wb{A}_{i,j} \times_Y \wb{B}_{i,j}]$ or
(2)  correspondence $\Gamma''$ such that $\dim(p_1|\Gamma''|) < \dim(X) -r$.

We now consider the action of such products on $H^*_c(X, \ql)$: By Lemma \ref{lem:compos}
and the fact that all our correspondences are supported in $X \times_Y X$ it follows
from the projection formula that the image of any product containing a factor of type
(1) must be contained in the subspace of $H^*_c(X, \ql)$ of the form 
$H^*_c(Y, \ql) \cdot \alpha$ where $\alpha$ denotes the class of some algebraic cycle
on $X$ in $H^*(X, \ql)$. It follows from Lemma \ref{lem:div}
that the eigenvalues of the geometric Frobenius
acting on the image of any product with a factor of type (2) are all divisible
by $q^{r+1}$.

Let 
$\alpha_{i,k} \in H^{2i}(X, \ql)$, $i =0,1, \dots, r$, $k = 1, \dots, b_{2i}$,
be algebraic cohomology classes such that their restrictions
 to the generic fibre $Z$ form a basis for
$\oplus_{i=0}^{r} H^{2i}(Z, \ql)$; these exist because of the assumptions on the
Chow groups of $Z$. By Lemma \ref{lem:inj} we have an injective
map
\[
\big ( \bigoplus_{i,k} H^*_c(Y, \ql) \big ) \sr{\oplus f^* \cup \alpha_{i,k}}{\lr} H^*_c(X, \ql)
\]
and the image of any product as above of type (1) is contained in the image
of this map. The above discussion then shows that all the eigenvalues of geometric Frobenius
acting on the cokernel of this map are  divisible by $q^{r+1}$.

We conclude the proof as in that of \cite[Theorem 1.1]{nf-rajan}: without loss
of generality we may assume that $Y$ has a unique rational point and then we
get the desired congruence by applying the Grothendieck-Lefschetz trace formula
to $X$.
\end{proof}

\begin{cor}
\label{cor:main}
Retain the hypotheses of Theorem \ref{thm:main}. Then  for all  $y \in Y(k)$,
all the eigenvalues of the geometric Frobenius acting on $H^i(Z_y, \ql)$
are divisible by $q^{r+1}$ for all $i > 2r$.

\end{cor}

\begin{proof}
This follows immediately from the proof of the theorem using the method 
of Esnault in \cite[Theorem 2.1]{esnault-app}; the idea is to use Artin's
vanishing theorem for the \'etale cohomology of affine varieties together
with the long exact sequence of cohomology with compact supports.
\end{proof}

\begin{rem}
All our arguments are essentially motivic in nature and can be easily adapted to
obtain analogous results for other cohomology theories; in particular the Hodge
theoretic analogues of Theorem \ref{thm:main} and Corollary \ref{cor:main}
are also true. See \cite{nf-rajan}, \cite{esnault-app} for the case $r=0$.
\end{rem}

\section{An application} \label{sec:examples}
We now give a family of examples to which our main theorem applies. The
proof of the proposition below is a variant of the proof of Otwinowska 
\cite[Corollaire 3]{Otwinowska}; the  novelty here is that we consider 
all hypersurfaces, not just the smooth ones.

In the following proposition all the data is assumed to be
defined over a universal domain $K$.
\begin{prop}
Let $B$ be a smooth projective variety with $Pic(B) \cong  \Z$ generated
by a very ample line bundle $L$ and such that the
cycle class maps $CH_i(B) \otimes \ql \to H_{2i}(B, \ql)$
are isomorphisms for $0 \leq i \leq 2 \dim(B)$. Let $d$ be such that for any section
$s \in H^0(B, {L}^{\otimes d})$, $Z(s)$, the zero scheme of $s$, is covered
by a family of projective spaces 
of dimension $k$ for some $k < \dim(B)/4$ and
 such that $L$ restricts to $\mc{O}(1)$ on each of them. 
Then for any $s$ as above, the cycle class maps
$CH_i(Z(s)) \otimes \ql \to H_{2i}(Z(s), \ql)$ for
$0 \leq i \leq k-1$ are isomorphisms. \label{prop:hyp}
\end{prop}

\begin{proof}
Let $A = Z(s)$ and first assume that $A$ is smooth. We recall the proof 
in this case, in a slightly modified form, from \cite{Otwinowska} since 
we need to extend it 
to allow for singularities: Let $G$ be the Grassmannian of $k$-planes
in $\P(H^0(B,L))$. The hypotheses, along with 
the use of  alterations \cite{dejong-alt1}, imply that we
may find a smooth projective variety $F$ with a morphism to $G$ such
that there is a dominant generically finite morphism $q:P \to A$,
with $p:P \to F$ the pullback of the universal bundle of $k$-planes on 
$G$ to $F$.

Let $h = c_1(\mc{O}_P(1)) \in CH^1(P)$. Then for any $i$ we have 
\begin{equation}
CH_{i}(P) = p^*CH_{i-k}(F) \oplus \big ( \bigoplus_{j=1}^k h^j\cdot p^* CH_{i -k +j}(F) \big) .
\end{equation}
For $Z \in CH_i(A)$  we may write
\[
q^*(Z) = Z_0 + h \cdot Z_1 + \cdots + h^k\cdot Z_k
\]
with $Z_j \in p^*CH_{i -k +j}(F)$. 
 
Let $H = c_1(L|_A) \in CH^1(A)$, so $h = q^*H$. Suppose
$i < k$ and $Z \in CH_i(A)$ is homologically trivial. Then in the above decomposition
$Z_0 = 0$ and all the $Z_j$'s are homologically trivial.
By the projection formula, $q_*(h\cdot Z_j) = H\cdot q_*(Z_j)$ which is zero
by the self-intersection formula, since, by the assumptions on $B$,
$i_*(q_*(Z_j)) = 0 \in CH_{i+j}(B)$. Thus $q_*(q^*(Z)) = 0$ which implies that
$Z = 0$ since we are working with $\Q$ coefficeients. The isomorphisms in
the theorem now follow becuase of the assumptions on
 $B$ and the  Lefschetz theorems.

Now we assume that $A$ is reduced and irreducible but not necessarily smooth.
We may find $F$ as above, but we do not in general have a map 
$q^*:CH(A) \to CH(P)$ such that $q^*$ sends homologically trivial cycles
to homologically trivial cycles and $q_*q^*$ is multiplication by the degree
of $q$. Nevertheless, $q_*:CH(P) \to CH(A)$ is surjective so for any $Z$
in $CH_i(A)$ with $i < k$ we may find $Z_j \in p^*CH_{i-k +j}(P)$ such that 
\[
Z = q_*(Z_1)\cdot H + \cdots + q_*(Z_k)\cdot H^k .
\]
By the self-intersection formula as above, the image of 
$CH_{i+j}(A) \stackrel{\cdot H^j}{\longrightarrow} CH_i(A)$ factors through 
$CH_{i+j}(A)/\equiv_{hom}$. By the assumptions on $B$ and Artin's theorem
on the cohomology of affine varieties \cite[Corollaire 3.2, p.160]{SGA4III}
it follows that the maps $H_{2j}(A,\ql) \to H_{2j}(B, \ql)$ are isomorphisms
for  $0 \leq j < 2k$.


Now suppose that  $A'$ is a smooth hypersurface in
$B$ specialising to  $A$, \emph{i.e}.~there is a commutative diagram
\[
\xymatrix{
\mc{A} \ar[rr] \ar[dr] & & \P(H^0(B, L^d)) \ar[dl]\times_K Spec(R) \\
& Spec(R) & \\
}
\]
where $R$ is a dvr with residue field  $K$ , $\mc{A} \to Spec(R)$ is a flat proper scheme
with geometric generic fibre isomorphic to $A'$ and special fibre isomorphic
to $A$ and the horizontal arrow is a closed embedding; such a diagram exists since
$L$ is very ample.

Using specialisation and statement for the smooth case, it follows that  the maps
 $CH_j(A)/\equiv_{hom}\otimes \ql \to H_{2j}(A, \ql)$ are surjective, hence
isomorphisms, for all $0 \leq j < 2k$, which implies that the specialisation maps 
$CH_j(A')/\equiv_{hom} \to CH_j(A)/\equiv_{hom}$ are also
isomorphisms;

By the above discussion and the fact that $c_1(L|_{A'})$ specializes
to $c_1(L|_A)$, we deduce  that for $0 \leq j < k$,
the specialisation maps $CH_j(A') \to CH_j(A)$ are surjective,  hence must be 
isomoprhisms since $\dim(H_{2j}(A',\ql)) = \dim(H_{2j}(A, \ql))$.

For the general case, let $A = \cup_{r=1}^m A_r$, where all the $A_r$'s are irreducible
divisors. We prove the proposition by induction on $m$, the case $m = 1$ being 
already done. Suppose $m=2$. We may find smooth divisors $A_1'$ and $A_2'$
specialising (as above) to $A_1$ and $A_2$ respectively such that $A_1'$ and 
$A_2'$ intersect transversally. By Artin's theorem  as before it follows that the map 
$H_{2j}(A_1' \cap A_2', \ql) \to H_{2j}(A_r', \ql)$
is surjective for $r = 1,2$, $0 \leq j < k$, so
that $CH_{j}(A_1' \cap A_2') \to CH_j(A_r')$
is also surjective. 
By specialisation (applied to 
 $A_1' \cap A_2'$ and  $A_1 \cap A_2$)
 it follows that $CH_j(A_1 \cap A_2) \to CH_j(A_r)$
is also surjective. This implies that $CH_j(A_r) \to CH_j(A_1 \cup A_2)$ is also
surjective. Since the analogous map on homology is an isomorphism by Artin's
theorem, it follows that
the above map on Chow groups is an isomorphism, concluding the proof in this case. 
The case $m > 2$ follows easily from the $m=2$ case.
\end{proof}

We will apply the above proposition in the case of certain homogenous spaces of
Picard rank one embedded into projective space via the very ample generator.
For that we shall need the following:
\begin{lem}
Let $k$ be an algebraically closed field and let $G_n$ be one of the simple
algebraic groups $SL_{n+1}, Sp_{2n}, SO_{2n}, SO_{2n+1}$ over $k$ and $P_n$ a
maximal parabolic subgroup of $G_n$ except for
the case where $G_n$ is $Sp_{2n}$ and $P_n$ is the stabilizer of a maxmimal isotropic
subspace.  Let $X_n = G_n/P_n$ and $L_n$ the (very) ample generator of $Pic(X_n)$.
Then  $X_n$ is coverered by a family of projective spaces of dimension
$m \geq (n-1)/2$
such that the restriction of $L_n$ to each of these subspaces is $\mc{O}(1)$.
\label{lem:linsubsp}
\end{lem}

\begin{proof}
This is well known and follows from the description of the $X_n$'s as spaces
of (isotropic) linear subspaces of the defining representation of $G_n$.
Note that in the case of the maximal isotropic subspaces for the orthogonal
groups, for the embedding induced by the Pl\"ucker embedding the restriction
of $\mc{O}(1)$ is not the generator so this case requires a slight
additional argument.
\end{proof}

We can now  prove our asymptotic generalisation of the theorem of Ax--Katz.

\begin{thm}
\label{thm:asaxkatz}
Let $k$ be  finite field with $q$ elements, let
$G_n$ be a $k$-form of one of the simple algebraic groups $SL_{n+1}$, $Sp_{2n}$, $SO_{2n}$
or $SO_{2n+1}$ and let
$P_n$ be a maximal parabolic subgroup of $G_n$ also defined over $k$ except for
the case where $G_n$ is $Sp_{2n}$ and $P_n$ is the stabilizer of a maxmimal isotropic
subspace. 
Let $X_n = G_n/P_n$ and $L_n$ the ample generator of $Pic(X_n)$. 
For any tuple  $D := (d_1,d_2, \ldots, d_t)$  of positive integers and any
positive integer $r$, there exists an integer $N(D,r)$ such that for 
any $n > N(D,r)$ and any $s_i \in H^0(X_n, L_n^{\otimes d_i})$, $i=1,2, \ldots,t$
\[
|Z(s_1,s_2,\ldots,s_t)(k)| \equiv \sum_{j=0}^{r-1}\dim(b_2(X_n))q^j \mod q^r
\]
where $Z(s_1,s_2,\ldots,s_t)$ is the subscheme of $X_n$ given by the vanishing of all
the $s_i$'s.
\end{thm}

\begin{proof}
Let $K \supset k$ be a universal domain. By applying Lemma \ref{lem:linsubsp}
to $X_n$ it follows that for fixed
positive integers $d,s$,
for any  $s \in H^0(X_{n,K}, L_{n,K}^{\otimes d})$, $Z(s)$ is
covered by linear spaces of dimension $s$ if $n >>0$, the reason being that
this is true for $X_n = \P^n$ (see e.g. \cite{paranjape}). 
By putting $A = X_n $ in Proposition \ref{prop:hyp} and using the fact that
for any projective homogenous space $X$ of a linear algebraic group the cycle
class map $CH(X)\te \ql \to H(X, \ql)$ is always an isomorphism,
we conclude that the cycle class map $CH_j(Z_s) \otimes \ql \to H_{2j}(Z(s), \ql)$
is an isomorphism for $j \leq s-1$. Moreover, $H_{2j+1}(Z(s), \ql)$ vanishes
for $j < k$ by the Lefschetz hyperplane theorem. 

By taking $s=r+1$, the previous paragraph shows that we may  apply Theorem 
\ref{thm:main} to conclude the proof for $t = 1$. Suppose $t > 1$: by the $t=1$ case
we may find an integer $N(D,r)$ so that for all $n> N(D,r)$ the desired
congruence holds for all the subschemes $\cup_{j \in T'} Z(s_j)$ where
$\emptyset \neq T' \subset \{1,2,\cdots,t\}$. The theorem then follows by
the Inclusion--Exclusion principle.
\end{proof}

\begin{ex} \label{ex:E8}
We consider a particular homogenous space for the exceptional group $G =E_8$;
our notation below is that of Bourbaki \cite[p. 268]{bourb-root}. Let $B$ be
a Borel subgroup of $G$ and let $P\supset B$ be
the parabolic subgroup of $G$ corresponding to (removing) the root $\alpha_1$; this 
is the stabilizer (upto scalars) of the highest weight vector $v$ in the representation
of $G$ with fundamental weight $\varpi_i$. Let $X = G/P$; using the data in
\emph{op.~cit}.~one easily sees that $\dim(X) = 78$. Let $P'\supset B$ be
the parabolic subgroup of $G$ corresponding to removing the root $\alpha_2$.
By looking at the Dynkin diagram of $E_8$ one sees that the semi-simple part 
$G'$ of a Levi subgroup of $P'$ is isomorphic to $SL_7$ and one
checks using the information in \emph{loc.~cit.} that the action of
$G'$ on $v$ generates the standard representation of $SL_7$. We thus see that
$X$ is covered by a family of $6$-dimensional projective spaces such that
the ample generator $L$ of $Pic(X)$ restricts to $\mc{O}(1)$.

By proceeding as  in the above theorem but using the explicit bound above
we get (where we now work over a finite field $k$ with $q$ elements):
\begin{enumerate}
\item  If $s \in H^0(X, L)$ then 
\[
|Z(s)(k)| \equiv 1 + \sum_{i=1}^{4}b_{2i}(X)q^i \mod q^5
\] \\
\item If $s \in H^0(X, L^{\otimes 2})$ then
\[
|Z(s)(k)| \equiv 1 + q \mod q^2 .
\]
This uses the fact that a quadric in $\P^6$ is covered by linearly embedded $\P^2$'s.
\end{enumerate}
\end{ex}

\begin{rem} One can probably use the method of Paranjape  \cite{paranjape}
or Esnault--Levine--Viehweg \cite{ELV} to prove a version of Proposition
\ref{prop:hyp} for all low degree intersections (instead of just hypersurfaces).
We have not done so for the sake of brevity since the proof is not likely to 
require any new ideas.
\end{rem}



\end{document}